\documentclass[12pt]{article}
\usepackage{amsmath,amsfonts,amsthm,amssymb,eucal}
\newcommand{\Qp}{\mathbb Q_p}
\newcommand{\DA}{D^\alpha}
\newcommand{\DAN}{D^\alpha_N}
\newcommand{\DD}{\mathcal D(\Qp)}
\newcommand{\DDD}{\mathcal D'(\Qp)}
\newcommand{\LL}{L^1(\mathbb Q_p)}
\newcommand{\FF}{\mathcal F}
\numberwithin{equation}{section}
\begin{document}
\newtheorem{lem}{Lemma}
\newtheorem*{teo}{Theorem}
\newtheorem{prop}{Proposition}
\newtheorem*{defin}{Definition}
\pagestyle{plain}
\title{$p$-Adic Analogue of the Porous Medium Equation}
\author{ \textbf{Andrei Yu. Khrennikov}\\
\footnotesize International Center for Mathematical Modeling\\
\footnotesize in Physics and Cognitive Sciences,\\
\footnotesize Mathematical Institute, Linnaeus University,\\
\footnotesize V\"axj\"o, SE-35195 Sweden,\\
\footnotesize E-mail: andrei.khrennikov@lnu.se
\and
\textbf{Anatoly N. Kochubei}\\
\footnotesize Institute of Mathematics,\\
\footnotesize National Academy of Sciences of Ukraine,\\
\footnotesize Tereshchenkivska 3, Kiev, 01004 Ukraine,\\
\footnotesize E-mail: kochubei@imath.kiev.ua }

\date{}
\maketitle

\bigskip
\begin{abstract}
We consider a nonlinear evolution equation for complex-valued functions of a real positive time variable and a $p$-adic spatial variable.
This equation is a non-Archimedean counterpart of the fractional porous medium equation. Developing, as a tool, an $L^1$-theory of Vladimirov's
$p$-adic fractional differentiation operator, we prove m-accretivity of the appropriate nonlinear operator, thus obtaining the existence and
uniqueness of a mild solution. We give also an example of an explicit solution of the $p$-adic porous medium equation.
\end{abstract}
\vspace{2cm}
{\bf Key words: }\ $p$-adic numbers; Vladimirov's $p$-adic fractional differentiation operator; $p$-adic porous medium equation;
mild solution of the Cauchy problem

\medskip
{\bf MSC 2010}. Primary: 35K55. Secondary: 11S80; 35R11.

\newpage

\section{Introduction}

The theory of linear partial pseudo-differential equations for complex-valued functions over non-Archimedean fields is a well-established branch of mathematical analysis. By this time, there is a description of various equations whose properties resemble those of classical equations of mathematical physics, there are constructions of fundamental solutions, information on spectral properties of related operators. For equations of evolution type, there are results on initial value problems etc. See \cite{AKS, AKS11,KS10,K2001,K2008,RZ,VVZ,Z2003,Z2008,Z2014} and references therein.

Meanwhile very little is known about nonlinear $p$-adic equations. We can mention only some semilinear evolution equations solved using $p$-adic wavelets \cite{AKS} and a kind of equations of reaction-diffusion type studied in \cite{Z2016}.

In this paper, we consider a $p$-adic analog of one of the most important classical nonlinear equations, the porous medium equation (see \cite{Va}), that is the equation
\begin{equation}
\label{1.1}
\frac{\partial u}{\partial t}+\DA (\varphi (u))=0,\quad u=u(t,x),\quad t>0,x\in \Qp,
\end{equation}
where $\Qp$ is the field of $p$-adic numbers, $\DA$, $\alpha >0$, is Vladimirov's fractional differentiation operator (see the deinitions and other preliminaries in Section 2 below), $\varphi$ is a strictly monotone increasing continuous real function, $|\varphi (s)|\le C|s|^m$ for $s\in \mathbb R$ ($C>0$, $m\ge 1$). A typical example of the latter is $\varphi (u)=u|u|^{m-1}$, $m>1$. We see Eq. (\ref{1.1}) as the simplest model example of an equation of this kind for the non-Archimedean situation. Therefore in order to understand specific features of this case, we confine ourselves to the simplest pseudo-differential operator on $\Qp$ and the simplest kind of nonlinearity. On the other hand, this setting has some common features with recent work on fractional porous medium equation on $\mathbb R^n$ \cite{PQRV}. Another motivation is the $p$-adic model of a porous medium proposed in \cite{KOC1,KOC2}.

Our strategy for studying Eq. (\ref{1.1}) is as follows. There exists an abstract theory of the equations
\begin{equation}
\label{1.2}
\frac{\partial u}{\partial t}+A (\varphi (u))=0.
\end{equation}
developed by Crandall and Pierre \cite{CP} and based on the theory of stationary equations
\begin{equation}
\label{1.3}
u+A\varphi (u)=f
\end{equation}
developed by Br\'ezis and Strauss \cite{BS}. In Eq. (\ref{1.2}) and (\ref{1.3}), $A$ is a linear $m$-accretive operator in $L^1(\Omega )$ where $\Omega$ is a $\sigma$-finite measure space. Under some natural assumptions, the nonlinear operator$A\varphi =A\circ \varphi$ is accretive and admits an $m$-accretive extension $A_\varphi$, the generator of a contraction semigroup of nonlinear operators. This result gives information on a kind of generalized solvability of Eq. (\ref{1.2}), though the available description of $A_\varphi$ is not quite explicit.

In order to use this method for Eq. (\ref{1.1}), we need an $L^1$-theory of the Vladimirov operator $\DA$. This is a subject of independent interest, and we treat it in Section 3.

In the classical situation where $\Omega =\mathbb R^n$, $A$ is the Laplacian, there are stronger results based on the study of Eq. (\ref{1.3}) (see \cite{BBC,C}) showing that $A\varphi$ is $m$-accretive itself. This employs some delicate tools of local analysis of solutions, such as imbedding theorems for Marcinkiewicz and Sobolev spaces in bounded domains.

For our $p$-adic situation, we prove (section 5) a little weaker result, namely the $m$-accretivity of the closure of the operator $A\varphi$. Our tool is the $L^1$-theory of the Vladimirov type operator on a $p$-adic ball (Section 4).

Finally, in Section 6, we give an example of an explicit solution of Eq. (\ref{1.1}) resembling the ``Quadratic Pressure Solution'' of the porous medium equation on $\mathbb R^n$ \cite{A}.

\medskip
\section{Preliminaries}

\medskip
{\it 2.1. p-Adic numbers and the Vladimirov operator} \cite{K2001,VVZ}.

Let $p$ be a prime
number. The field of $p$-adic numbers is the completion $\mathbb Q_p$ of the field $\mathbb Q$
of rational numbers, with respect to the absolute value $|x|_p$
defined by setting $|0|_p=0$,
$$
|x|_p=p^{-\nu }\ \mbox{if }x=p^\nu \frac{m}n,
$$
where $\nu ,m,n\in \mathbb Z$, and $m,n$ are prime to $p$. $\Qp$ is a locally compact topological field.

Note that by Ostrowski's theorem there are no absolute values on $\mathbb Q$, which are not equivalent to the ``Euclidean'' one,
or one of $|\cdot |_p$.

The absolute value $|x|_p$, $x\in \mathbb Q_p$, has the following properties:
\begin{gather*}
|x|_p=0\ \mbox{if and only if }x=0;\\
|xy|_p=|x|_p\cdot |y|_p;\\
|x+y|_p\le \max (|x|_p,|y|_p).
\end{gather*}

The latter property called the ultra-metric inequality (or the non-Archi\-me\-dean property) implies the total disconnectedness of $\Qp$ in the topology
determined by the metric $|x-y|_p$, as well as many unusual geometric properties. Note also the following consequence of the ultra-metric inequality:
\begin{equation*}
|x+y|_p=\max (|x|_p,|y|_p)\quad \mbox{if }|x|_p\ne |y|_p.
\end{equation*}

The absolute value $|x|_p$ takes the discrete set of non-zero
values $p^N$, $N\in \mathbb Z$. If $|x|_p=p^N$, then $x$ admits a
(unique) canonical representation \index{Canonical
representation}
\begin{equation}
\label{2.1}
x=p^{-N}\left( x_0+x_1p+x_2p^2+\cdots \right) ,
\end{equation}
where $x_0,x_1,x_2,\ldots \in \{ 0,1,\ldots ,p-1\}$, $x_0\ne 0$.
The series converges in the topology of $\mathbb Q_p$. For
example,
$$
-1=(p-1)+(p-1)p+(p-1)p^2+\cdots ,\quad |-1|_p=1.
$$
We denote $\mathbb Z_p=\{ x\in \Qp:\ |x|_p\le 1\}$. $\mathbb Z_p$, as well as all balls in $\Qp$, is simultaneously open and closed.

Proceeding from the canonical representation (\ref{2.1}) of an element $x\in
\mathbb Q_p$, we define the fractional part of $x$ as the rational number
$$
\{ x\}_p=\begin{cases}
0,& \text{if $N\le 0$ or $x=0$};\\
p^{-N}\left( x_0+x_1p+\cdots +x_{N-1}p^{N-1}\right) ,& \text{if
$N>0$}.\end{cases}
$$
The function $\chi (x)=\exp (2\pi i\{ x\}_p)$ is an additive
character of the field $\mathbb Q_p$, that is a character of its additive group. It is clear
that $\chi (x)=1$ if $|x|_p\le 1$. Denote by $dx$ the Haar measure on the
additive group of $\Qp $ normalized by the equality $\int_{\mathbb Z_p}dx=1$.

The above additive group is self-dual, so that
the Fourier transform of a complex-valued function $f\in L^1(\Qp)$ is again a function on $\Qp$ defined as
$$
(\mathcal Ff)(\xi )=\int\limits_{\Qp}\chi (x\xi )f(x)\,dx.
$$
If $\mathcal Ff\in L_1(\Qp)$, then we have the inversion formula
$$
f(x)=\int\limits_{\Qp}\chi (-x\xi )\widetilde{f}(\xi )\,d\xi .
$$
It is possible to extend $\mathcal F$ from $L_1(\Qp)\cap L_2(\Qp)$ to a unitary operator on $L_2(\Qp)$, so that the Plancherel identity holds in this case. In a way, Fourier analysis on $\Qp$ is close to the classical harmonic analysis (see \cite{Fei} for an interesting discussion of various stages in the development of harmonic analysis).

In order to define distributions on $\Qp$, we need a class of test functions. A function $f:\ \Qp\to \mathbb C$ is called locally constant if
there exists such an integer $l\ge 0$ that for any $x\in \Qp$
$$
f(x+x')=f(x)\quad \mbox{if }\|x'\|\le p^{-l}.
$$
The smallest number $l$ with this property is called the exponent of local constancy of the function $f$.

Typical examples of locally constant functions are additive characters, and also cutoff functions like
$$
\Omega (x)=\begin{cases}
1,& \text{if $\|x\|\le 1$};\\
0,& \text{if $\|x\|>1$}.\end{cases}
$$
In particular, $\Omega$ is continuous, which is an expression of the non-Archimedean properties of $\Qp$.

Denote by $\mathcal D(\Qp)$ the vector space of all locally constant functions with compact supports. Note that $\DD$ is dense in $L^q(\Qp)$ for each $q\in [1,\infty )$. In order to furnish $\DD$ with a topology, consider first the subspace $D_N^l\subset \DD$ consisting of functions with supports in a ball
$$
B_N=\{ x\in \Qp:\ |x|_p\le p^N\},\quad N\in \mathbb Z,
$$
and the exponents of local constancy $\le l$. This space is finite-dimensional and possesses a natural direct product topology. Then the topology in $\DD$ is defined as the double inductive limit topology, so that
$$
\DD=\varinjlim\limits_{N\to \infty}\varinjlim\limits_{l\to \infty}D_N^l.
$$

If $V\subset \Qp$ is an open set, the space $\mathcal D(V)$ of test functions on $V$ is defined as a subspace of $\DD$ consisting of functions with supports in $V$.

The space $\DDD$ of Bruhat-Schwartz distributions on $\Qp$ is defined as a strong conjugate space to $\DD$.

In contrast to the classical situation, the Fourier transform is a linear automorphism of the space $\DD$. By duality, $\FF$ is extended to a linear automorphism of $\DDD$. There exists a detailed theory of convolutions and direct products of distributions on $\Qp$ closely connected with the theory of their Fourier transforms; see \cite{AKS,K2001,VVZ}.

The Vladimirov operator $\DA$, $\alpha >0$, of fractional differentiation, is defined first as a pseudo-differential operator with the symbol $|\xi_p^\alpha$:
\begin{equation}
\label{2.2}
(\DA u)(x)=\FF^{-1}_{\xi \to x}\left[ |\xi |_p^{\alpha }\FF_{y\to \xi }u\right] ,\quad u\in \DD,
\end{equation}
where we show arguments of functions and their direct/inverse Fourier transforms. There is also a hypersingular integral representation giving the same result on $\DD$ but making sense on much wider classes of functions (for example, bounded locally constant functions):
\begin{equation}
\label{2.3}
\left( \DA u\right) (x)=\frac{1-p^\alpha }{1-p^{-\alpha -1}}\int\limits_{\Qp}|y|_p^{-\alpha -1}[u(x-y)-u(x)]\,dy.
\end{equation}

The Cauchy problem for the heat-like equation
$$
\frac{\partial u}{\partial t}+\DA u=0,\quad u(0,x)=\psi (x),\quad x\in\Qp,t>0,
$$
possesses many properties resembling classical parabolic equations. If $\psi$ is regular enough, for example, $\psi \in \DD$, then a classical solution is given by the formula
$$
u(t,x)=\int\limits_{\Qp}Z(t,x-\xi )\psi (\xi )\,d\xi
$$
where $Z$ is, for each $t$, a probability density and
\begin{equation}
\label{2.4}
Z(t_1+t_2,x)=\int\limits_{\Qp}Z(t_1,x-y)Z(t_2,y)\,dy,\quad t_1,t_2>0,\ x\in \Qp.
\end{equation}
Explicitly,
\begin{equation}
\label{2.5}
Z(t,x)=\sum\limits_{k=-\infty}^\infty p^kc_k(t)\Delta_{-k}(x)
\end{equation}
where $\Delta_l(x)$ is the indicator function of the ball $B_l$,
\begin{equation}
\label{2.6}
c_k(t)=\exp \left( -p^{k\alpha }t\right)-\exp \left( -p^{(k+1)\alpha }t\right).
\end{equation}

Another expression for $Z(t,x)$, valid for $x\ne 0$, is
\begin{equation}
\label{2.7}
Z(t,x)=\sum\limits_{m=1}^\infty \frac{(-1)^m}{m!}\cdot \frac{1-
p^{\alpha m}}{1-p^{-\alpha m-1}}t^m|x|_p^{-\alpha m-1}.
\end{equation}

The ``heat kernel'' $Z$ satisfies the estimate
\begin{equation}
\label{2.8}
0<Z(t,x)\le Ct(t^{1/\alpha }+|x|_p)^{-\alpha -1},\quad t>0,x\in \Qp.
\end{equation}
Here and below the letter $C$ denotes various positive constants.

\bigskip
{\it 2.2. A Heat-Like Equation on a p-Adic Ball} \cite{K2001}

Let us consider the Cauchy problem
\begin{gather}
\frac{\partial u(t,x)}{\partial t}+\left( D_N^\alpha u\right) (t,x)-\lambda
u(t,x)=0,\quad x\in B_N,\ t>0;\label{2.9}\\
u(0,x)=\psi (x),\quad x\in B_N,\label{2.10}
\end{gather}
where $N\in \mathbb Z$, $\psi\in \mathcal D(B_N)$, $\lambda =\dfrac{p-1}{p^{\alpha +1}-1}p^{\alpha (1-N)}$, the operator $D_N^\alpha$ is defined by restricting $\DA$ to functions $u_N$ supported in $B_N$ and considering the resulting function $\DA u_N$ only on $B_N$. Here and below we often identify a function on $B_N$ with its extension by zero onto $\Qp$. Note that $D_N^\alpha$ defines a positive definite operator on $L^2(B_N)$, $\lambda$ is its smallest eigenvalue.

The probabilistic meaning of the problem (\ref{2.9})-(\ref{2.10}) is discussed in \cite{K2001}. Here we will write only its solution
\begin{equation*}
u(x,t)=\int \limits _{B_N}Z_N(t,x-y)\psi (y)\,dy,\quad t>0,x\in B_N,
\end{equation*}
where
\begin{gather}
Z_N(t,x)=e^{\lambda t}Z(t,x)+c(t),\quad x\in B_N,\label{2.11}\\
c(t)=p^{-N}-p^{-N}(1-p^{-1})e^{\lambda t}\sum \limits _{n=0}^\infty
\frac{(-1)^n}{n!}t^n\frac{p^{-N\alpha n}}{1-p^{-\alpha n-1}},\label{2.12}
\end{gather}
The kernel $Z_N$ is a transition density of a Markov process on $B_N$.

\bigskip
{\it 2.3. Nonlinear Semigroups and the Abstract Porous Medium Equation} \cite{Ba,BS,CP}.

Let $\Omega$ be a $\sigma$-finite measure space, $A:\ D(A)\subset L^1(\Omega)\to L^1(\Omega)$ be a linear densely defined operator, such that $A$ is $m$-accretive, that is a generator of a strongly continuous contraction semigroup $e^{-tA}$. It is assumed, in addition, that
\begin{equation}
\label{2.13}
0\le f\le 1\Longrightarrow 0\le e^{-tA}f\le 1.
\end{equation}

Let $\varphi :\ \mathbb R\to \mathbb R$ be a continuous strictly increasing function, $\varphi (0)=0$. Consider the nonlinear operator $A\varphi$ with the domain
$$
D(A\varphi )=\{ u\in L^1(\Omega ):\ \varphi (u)\in D(A)\}.
$$
It is proved \cite{CP} that $(A\varphi )(u)=A(\varphi (u))$ is an accretive nonlinear operator in $L^1(\Omega )$, and for any $\varepsilon >0$, $(\varepsilon I+A)\varphi$ is $m$-accretive for each $\lambda >0$, $(I+\lambda A\varphi)^{-1}$ is an order-preserving non-expansive mapping. However it may happen that $A\varphi$ is not $m$-accretive.

On the other hand, $A\varphi$ has a $m$-accretive extension $A_\varphi$ in $L^1(\Omega )$, which extends $A\varphi$ in such a way that for any $\lambda >0$, $f\in L^1(\Omega )$,
$$
(I+\lambda A_\varphi )^{-1}f=\lim\limits_{\varepsilon \downarrow 0}(I+\lambda (\varepsilon I+A)\varphi )^{-1}f.
$$

This result means a kind of generalized solvability of the Cauchy problem for the equation
\begin{equation}
\label{2.14}
\frac{\partial u}{\partial t}+A (\varphi (u))=0.
\end{equation}
The $m$-accretivity of a nonlinear operator $B$ means that it determines a nonlinear strongly continuous nonexpansive semigroup corresponding in a usual way to the Cauchy problem.

Another option is to consider the closure of the operator $A\varphi$, that is the closure of its graph, which is a single-valued operator, since $A$ is closed and $\varphi$ is continuous. The closure $\overline{A\varphi}$ is obviously accretive. It is $m$-accretive, if the range of $I+A\varphi$ is dense in $L^1(\Omega )$ or, in other words, the equation
$$
u+A\varphi (u)=f
$$
is solvable for a dense subset of functions $f\in L^1(\Omega )$.

Equivalently, setting $\beta =\varphi^{-1}$ (the function inverse to $\varphi$), we have to study the equation
\begin{equation}
\label{2.15}
Av+\beta (v)=f.
\end{equation}
Together with Eq. (\ref{2.15}), one considers (see \cite{BS}) the regularized equation
\begin{equation}
\label{2.16}
\varepsilon v_\varepsilon +Av_\varepsilon +\beta (v_\varepsilon )=f,\quad \varepsilon >0,
\end{equation}
possessing a unique solution $v_\varepsilon$, such that $w_\varepsilon=f-Av_\varepsilon$ satisfies the inequality
\begin{equation}
\label{2.17}
\|w_\varepsilon\|_{L^1(\Omega )}\le \|f\|_{L^1(\Omega )}.
\end{equation}

Moreover, if $\hat{v}_\varepsilon$ and $\hat{w}_\varepsilon$ correspond to Eq. (\ref{2.16}) with a different right-hand side $\hat{f}$, then
\begin{equation}
\label{2.18}
\|w_\varepsilon -\hat{w}_\varepsilon\|_{L^1(\Omega )}\le \|f-\hat{f}\|_{L^1(\Omega )}.
\end{equation}
In addition, if $f\in L^1(\Omega)\cap L^\infty (\Omega)$, then
\begin{equation}
\label{2.19}
\|f-(A+\varepsilon)v_\varepsilon\|_{L^\infty (\Omega )}= \|\beta (v_\varepsilon )\|_{L^\infty (\Omega )}\le \|f\|_{L^\infty(\Omega )}
\end{equation}
(see Proposition 4 in \cite{BS}). These inequalities will be used in the study of Eq. (\ref{1.1}) in Section 5.

\section{The Vladimirov Operator in $\LL$}

{\it 3.1. The Heat-Like Equation and the Corresponding Semigroup of Operators}.

Using the fundamental solution $Z$ described in Section 2.1, we define the operator family
$$
(S(t)\psi )(x)=\int\limits_{\Qp}Z(t,x-\xi )\psi (\xi )\,d\xi,\quad \psi\in \LL,
$$
$t>0$. It follows from (\ref{2.4}), (\ref{2.8}) and the Young inequality that $S$ is a contraction semigroup in $\LL$.

\medskip
\begin{prop}
$S(t)$ has the $C_0$-property.
\end{prop}

\medskip
{\it Proof}. Since the space $\mathcal D(\Qp)$ of Bruhat-Schwartz test functions is dense in $\LL$ \cite{VVZ}, it follows from the equality
\begin{equation}
\label{3.1}
\int\limits_{\Qp}Z(t,x)\,dx=1
\end{equation}
that it is sufficient to prove the limit relation
\begin{equation}
\label{3.2}
I_t=\int\limits_{\Qp}dx \int\limits_{\Qp}Z(t,x-\xi)|\psi (\xi)-\psi (x)|\,d\xi \longrightarrow 0,
\end{equation}
as $t\to 0$, for any $\psi \in \mathcal D(\Qp)$.

Suppose that $\psi (\xi)-\psi (x)=0$ whenever $|x-\xi|_p\le p^m$, and $\psi (x)=0$, if $|x|_p>p^N$, $m<N$. Let us write $I_t=I_t^{(1)}+I_t^{(2)}$ where $I_t^{(1)}$ corresponds to the integration in (\ref{3.2}) in $x,|x|_p\le p^N$, while $I_t^{(2)}=I_t-I_t^{(1)}$. Then
$$
I_t^{(1)}\le Ct\int\limits_{|x|_p\le p^N}dx\int\limits_{|x-\xi|_p>p^m}\left( t^{1/\alpha}+|x-\xi|_p\right)^{-\alpha -1}d\xi \\
\le Ct\int\limits_{|x|_p\le p^N}dx\int\limits_{|z|_p>p^m}|z|_p^{-\alpha -1}dz\to 0,
$$
as $t\to 0$.

Next, in $I_t^{(2)}$ the term $\psi (x)$ is absent, so that
$$
I_t^{(2)}=\int\limits_{|x|_p> p^N}dx\int\limits_{|\xi|_p\le p^N}Z(t,x-\xi)|\psi (\xi )|\,d\xi \le Ctp^N\int\limits_{|x|_p> p^N}\left( t^{1/\alpha}+|x|_p\right)^{-\alpha -1}dx
$$
because $|x-\xi|_p=|x|_p$ for $|x|_p>p^N$ and $|\xi|_p\le p^N$. We get
$$
I_t^{(2)}\le Ctp^N\int\limits_{|x|_p> p^N}|x|_p^{-\alpha -1}dx\longrightarrow 0,
$$
as $t\to 0$. $\qquad \blacksquare$

\medskip
\begin{defin}
We define the realization $A$ of $\DA$ in $\LL$ as the generator of the semigroup $S(t)$.
\end{defin}

Let $D(A)$ be the domain of the operator $A$. We will also use the following notation:
$$
B_{l,x_0}=\{ x\in \Qp:\ |x-x_0|_p\le p^l\}\quad (x_0\in \Qp),
$$
$\Delta_{l,x_0}$ is the indicator of the ball $B_{l,x_0}$, $\delta_l(x)=p^l\Delta_{-l}(x)$.

\medskip
\begin{prop}
If $u\in \mathcal D(\Qp)$, then $u\in D(A)$ and $Au=\DA u$ where the right-hand side is understood as usual in terms of the Fourier transform or the hypersingular integral representation.
\end{prop}

\medskip
{\it Proof}. An arbitrary function from $\mathcal D(\Qp)$ is a finite linear combination of the functions $\Delta_{l,x_0}$. Therefore it is sufficient to consider the case $u=\Delta_{l,x_0}$.

\medskip
\begin{lem}
There is the identity
\begin{equation}
\label{3.3}
\left( \Delta_{k,x_0}*\Delta_{l,x_1}\right) (x)=p^{k+l-\max (k,l)}\Delta_{\max (k,l)} (x-x_0-x_1).
\end{equation}
\end{lem}

\medskip
{\it Proof}. It is known that $\FF (\Delta_k)=\delta_k$ (\cite{V}, (12.1)), so that
$$
\FF(\Delta_{k,x_0})(x)=\int\limits_{|\xi -x_0|_p\le p^k}\chi (\xi x)\,d\xi=\chi (x_0x)\int\limits_{|\eta |_p\le p^k}\chi (\eta x)\,d\eta=\chi (x_0x)\delta_k(x).
$$

Next,
$$
\left( \FF \left( \Delta_{k,x_0}*\Delta_{l,x_1}\right) \right) (\xi)=\chi (x_0\xi )\chi (x_1\xi )\delta_k(\xi )\delta_l(\xi )=\chi ((x_0+x_1)\xi)p^{k+l}\Delta_{-\max (k,l)}(\xi ).
$$
Therefore
$$
\left( \Delta_{k,x_0}*\Delta_{l,x_1}\right) (x)=p^{k+l}\int\limits_{\Qp}\chi ((x_0+x_1-x)\xi)\Delta_{-\max (k,l)}(\xi )\,d\xi =p^{k+l}\delta_{-\max (k,l)}(x-x_0-x_1),
$$
which implies (\ref{3.3}). $\qquad \blacksquare$

\medskip
\begin{lem}
There is the identity
\begin{equation}
\label{3.4}
(Z(t,\cdot )*\Delta_{l,x_0})(x)=\exp (-ap^{-\alpha l}t)\Delta_l(x-x_0)+p^l\sum\limits_{k=l+1}^\infty c_{-k}(t)p^{-k}\Delta_k(x-x_0).
\end{equation}
\end{lem}

\medskip
{\it Proof}. By (\ref{2.5}), (\ref{2.6}) and (\ref{3.3}), we have
\begin{multline*}
(Z(t,\cdot )*\Delta_{l,x_0})(x)=\sum\limits_{k=-\infty}^\infty p^kc_k(t)(\Delta_{-k}*\Delta_{l,x_0})(x)=\sum\limits_{k=-\infty}^\infty c_k(t)p^{l-\max (-k,l)}\Delta_{\max (-k,l)}(x-x_0) \\
=\left[ \sum\limits_{k=-l}^\infty c_k(t)\right] \Delta_l(x-x_0)+\sum\limits_{k=-\infty}^{-l-1}c_k(t)p^{l+k}\Delta_{-k}(x-x_0)=\exp (-ap^{-\alpha l}t)\Delta_l(x-x_0)\\
+p^l\sum\limits_{k=-\infty}^{-l-1}c_k(t)p^k\Delta_{-k}(x-x_0),
\end{multline*}
which is equivalent to (\ref{3.4}). $\qquad \blacksquare$

\medskip
\begin{lem}
There is the identity
\begin{equation}
\label{3.5}
\left( \DA \Delta_{l,x_0}\right) (y)=\begin{cases}
p^l\frac{1-p^{-1}}{1-p^{-\alpha -1}}p^{-l(\alpha +1)}, & \text{ if $|y-x_0|_p\le p^l$};\\
p^l\Gamma_p(\alpha +1)|y-x_0|_p^{-\alpha -1}, & \text{if $|y-x_0|_p> p^l$}.
\end{cases}
\end{equation}
\end{lem}

\medskip
{\it Proof}. We have
$$
\left( \DA \Delta_{l,x_0}\right) (y)=\FF^{-1}\left( |\xi |_p^\alpha \chi (x_0\xi )\delta_l(\xi )\right) (y)=p^l\int\limits_{\Qp}\chi (-(y-x_0)\xi )|\xi |_p^\alpha \Delta_{-l}(\xi )\,d\xi ,
$$
and it remains to use the identities (12.38) and (12.39) from \cite{V}. $\qquad\blacksquare$

\medskip
\begin{lem}
The coefficients $c_{-k},k\ge l+1$, admit the representation
\begin{equation}
\label{3.6}
c_{-k}(t)=p^{-k\alpha }(p^\alpha -1)t+b_kt^2,\quad k\ge l+1,
\end{equation}
where $|b_k|\le Cp^{-2k\alpha}$, $k\ge l+1$.
\end{lem}

\medskip
{\it Proof}. In accordance with (\ref{2.6}),
$$
c_{-k}(t)=e^{-p^{-k\alpha} t}-e^{-p^{-k\alpha +\alpha}t}.
$$
We have $c_{-k}(0)=0$,
$$
c_{-k}'(t)=-p^{-k\alpha}e^{-p^{-k\alpha} t}+p^{-k\alpha +\alpha}e^{-p^{-k\alpha +\alpha}t},
$$
$$
c_{-k}''(t)=p^{-2k\alpha}e^{-p^{-k\alpha} t}-p^{-2k\alpha +2\alpha} e^{-p^{-k\alpha +\alpha}t},
$$
so that
$$
\left| c_{-k}''(t)\right| \le Cp^{-2k\alpha} \quad \text{for all} t\ge 0.
$$
where the constant $C>0$ does not depend on $k,t$.

By the Taylor formula,
$$
c_{-k}(t)=tc_{-k}'(0)+t^2b_k,\quad |b_k|\le Cp^{-2k\alpha},
$$
which implies (\ref{3.6}). $\qquad \blacksquare$

\medskip
{\it Proof of Proposition 2 (continued)}. We have to consider the expression ($u=\Delta_{l,x_0}$)
\begin{multline*}
\left( \frac1t[-S(t)u+u]-\DA u\right) (x)=\frac1t \left( -e^{-p^{\alpha l}t}+1\right) \Delta_l(x-x_0)\\
-(p^\alpha -1)p^l\sum\limits_{k=l+1}^\infty p^{-k(\alpha +1)}\Delta_k(x-x_0)-tp^l\sum\limits_{k=l+1}^\infty p^{-k}b_k\Delta_k(x-x_0)
-\left( \DA \Delta_{l,x_0}\right) (x).
\end{multline*}

We have
$$
\frac1t \left( e^{-p^{\alpha l}t}-1\right) =-p^{-\alpha l}+\omega (t),
$$
$\omega (t)\to 0$, as $t\to 0$. If $|x-x_0|_p\le p^l$, then by (\ref{3.5}),
\begin{multline*}
p^{-\alpha l}\Delta_l(x-x_0)-(p^\alpha -1)p^l\sum\limits_{k=l+1}^\infty p^{-k(\alpha +1)}\Delta_k(x-x_0)-\left( \DA \Delta_{l,x_0}\right) (x)\\
=p^{-\alpha l}-(p^\alpha -1)p^l\sum\limits_{k=l+1}^\infty p^{-k(\alpha +1)}-p^l\frac{1-p^{-1}}{1-p^{-\alpha -1}}p^{-l(\alpha +1)}=0,
\end{multline*}
which is checked by an elementary calculation.

If $|x-x_0|_p=p^m$, $m\ge l+1$, then
$$
\Delta_k(x-x_0)=\begin{cases}
0, & \text{if $k<m$},\\
1, & \text{if $k\ge m$},
\end{cases}
$$
so that
$$
(p^\alpha -1)p^l\sum\limits_{k=l+1}^\infty p^{-k(\alpha +1)}\Delta_k(x-x_0)=(p^\alpha -1)p^l\sum\limits_{k=m}^\infty p^{-k(\alpha +1)}
=-p^l\Gamma_p(\alpha +1)|x-x_0|_p^{-\alpha -1},
$$
and the sum of this expression with the one for $\left( \DA \Delta_{l,x_0}\right) (x)$ for $|x-x_0|_p=p^m$ equals zero.

Therefore
\begin{multline*}
\left\| \frac1t [-S(t)\Delta_{l,x_0}+\Delta_{l,x_0}]-\DA \Delta_{l,x_0}\right\|_{\LL}\le |\omega (t)|\int\limits_{\Qp}\Delta_l(x-x_0)\,dx\\
+tp^l \sum\limits_{k=l+1}^\infty p^{-k}b_k\int\limits_{\Qp}\Delta_k(x-x_0)\,dx=p^l|\omega (t)|+tp^l\sum\limits_{k=l+1}^\infty b_k\longrightarrow 0,
\end{multline*}
as $t\to 0$, as desired. $\qquad\blacksquare$

\bigskip
{\it 3.2. The Green function.}

Since the operator $A$ in $\LL$ is defined as the generator of the contraction semigroup $S(t)=e^{-tA}$, then by the Hille-Yosida theorem, we can find the resolvent $R_\mu (A)=(A+\mu I)^{-1}$, $\mu >0$, by the formula
\begin{equation}
\label{3.7}
R_\mu (A)\psi =-\int\limits_0^\infty e^{-\mu t}S(t)\psi \,dt,\quad \psi \in \LL.
\end{equation}

Using (\ref{2.5}) and (\ref{3.7}) we first write
$$
(S(t)\psi )(x)=\sum\limits_{k=-\infty}^\infty p^k\left( e^{-p^{-k\alpha} t}-e^{-p^{-k\alpha +\alpha}t}\right)\int\limits_{|\eta |_p\le p^{-k}}\psi (x-\eta )\,d\eta,
$$
and then obtain the equality
\begin{equation}
\label{3.8}
(R_\mu (A)\psi)(x)=(p^\alpha -1)\sum\limits_{k=-\infty}^\infty p^{k(\alpha +1)}\frac1{(\mu +p^{\alpha k})(\mu +p^{\alpha k+\alpha})}\int\limits_{|\eta |_p\le p^{-k}}\psi (x-\eta )\,d\eta.
\end{equation}

The convergence in the right-hand side of (\ref{3.8}) is obvious for $k\to -\infty$, since $\psi \in \LL$. The convergence, as $k\to \infty$, for almost all $x\in \Qp$, follows from the relation
$$
p^k\int\limits_{|\eta |_p\le p^{-k}}\psi (x-\eta )\,d\eta =p^k\int\limits_{|x-\xi|_p\le p^{-k}}[\psi (\xi )-\psi (x)]\,d\xi +\psi (x)\to \psi (x),\quad k\to \infty,
$$
due to the theorem about Lebesgue points (proved for general measure spaces in \cite{Fe}, Theorem 2.9.8).

We will consider below the case where $\alpha >1$, in which the resolvent is an integral operator with a kernel possessing some smoothness properties. Thus, from now on,
\begin{equation}
\label{3.9}
\alpha >1.
\end{equation}
It is proved (for a more general situation) in \cite{K2001} (Section 2.5.3) that in this case $R_\mu$ is a convolution operator with the continuous integral kernel $E_\mu (x-\xi )$, such that $E_\mu (x)\sim \operatorname{const}\cdot |x|_p^{-\alpha -1}$, $|x|_p\to \infty$ (note a misprint in the formula (2.25) of \cite{K2001}). The function $E_\mu$ is represented by the uniformly convergent series
\begin{equation}
\label{3.10}
E_\mu (x)=\sum\limits_{N=-\infty}^\infty e_\mu^{(N)}(x),
\end{equation}
\begin{equation}
\label{3.11}
e_\mu^{(N)}(x)=\int\limits_{|\xi|_p=p^N}\frac{\chi (-x\xi )}{|\xi|_p^\alpha +\mu}\,d\xi .
\end{equation}

To measure smoothness of $E_\mu$, we consider
$$
\| E_\mu (\cdot )-E_\mu (\cdot +h)\|_{\LL}, \quad |h|_p\le 1.
$$
Note that $E_\mu (x)$ depends actually on $|x|_p$ (see \cite{VVZ} or \cite{K2001}). If $|x|_p>1$, $|h|_p\le 1$, then $E_\mu (x)-E_\mu (x+h)=0$. Therefore
$$
\| E_\mu (\cdot )-E_\mu (\cdot +h)\|_{\LL}=\| E_\mu (\cdot )-E_\mu (\cdot +h)\|_{L^1(\mathbb Z_p)}.
$$
Denote this function of $h$ by $\Phi (h)$, $h\in \mathbb Z_p$. We have
$$
\Phi (h)=\int\limits_{\mathbb Z_p}dx \left| \int\limits_{\Qp}\chi (-x\xi )[1-\chi (-h\xi )]\frac{d\xi}{|\xi|_p^\alpha +\mu}\right| ,
$$
and by the dominated convergence theorem, if $\alpha >1$, then $\Phi (h)\to 0$, as $h\to 0$.

Note that classically \cite{Gu}, the characterization of smoothness is performed in $L^1$-theory in terms of the Besov space. In our situation, the existing theory of $p$-adic Besov spaces \cite{Ha,Ka,Ta} does not work so far.

If $f\in \LL$, $g=E_\mu *f$, then by the Young inequality
\begin{equation}
\label{3.12}
\| g(\cdot )-g(\cdot +h)\|_{\LL}\le \Phi (h)\|f\|_{\LL},
\end{equation}
so that the left-hand side of (\ref{3.12}) tends to zero, as $h\to 0$, uniformly with respect to $L^1$-bounded sets of functions $f$. This fact can be related to compactness criteria used in the theory of nonlinear equations.

\bigskip
{\it 3.3. Description of $A$ in the distribution sense.}

Let $u\in \LL$. Then $\DA u$ can be defined as a distribution from $\mathcal D'(\Qp)$, a convolution $u*f_{-\alpha}$, $f_{-\alpha}(x)=\dfrac{|x|_p^{-\alpha -1}}{\Gamma_p(-\alpha)}$,
\begin{equation}
\label{3.13}
\Gamma_p(z)=\frac{1-p^{z-1}}{1-p^{-z}}.
\end{equation}
By definition \cite{VVZ},
$$
\langle u*f_{-\alpha},\psi \rangle =\lim\limits_{k\to \infty}\langle u(x)\langle f_{-\alpha} (y),\Delta_k(y)\psi (x+y)\rangle_y\rangle_x,\quad \psi \in \mathcal D(\Qp).
$$
Here $f_{-\alpha}$ is defined by analytic continuation, which leads to the usual formula
\begin{multline*}
\langle f_{-\alpha} (y),\Delta_k(y)\psi (x+y)\rangle =\frac{1-p^\alpha}{1-p^{-\alpha -1}}\int\limits_{|y|_p\le p^k}|y|_p^{-\alpha -1}[\psi (x+y)-\psi (x)]\,dy\\
=\frac{1-p^\alpha}{1-p^{-\alpha -1}}\int\limits_{p^{-l}<|y|_p\le p^k}|y|_p^{-\alpha -1}[\psi (x+y)-\psi (x)]\,dy
\end{multline*}
where $l$ is the exponent of local constancy of the function $\psi$.

The last expression has a limit, as $k\to \infty$, uniform in $x\in \Qp$ and belonging to $L^\infty (\mathbb Q_p)$. Since $u\in \LL$, the operator $\DA$ is defined in the distribution sense:
$$
\langle \DA u,\psi \rangle =\langle u,\DA \psi \rangle ,\quad \psi \in \mathcal D(\Qp),
$$
where $\DA \psi \in L^\infty (\mathbb Q_p)$, $\langle u,v\rangle =\int\limits_{\Qp}u(x)\overline{v(x)}\,dx$.

\medskip
\begin{prop}
The operator $A$ defined as a semigroup generator has the domain $D(A)=\{ u\in \LL:\ \DA u\in \LL\}$ where $Au=\DA u$ (understood in the distribution sense).
\end{prop}

\medskip
{\it Proof}. Let
$$
u(x)=\int\limits_{\Qp}E_\mu (x-y)f(y)\,dy,\quad f\in \LL.
$$
Let us check that $(\DA +\mu I)u=f$ in the distribution sense. By the construction of $E_\mu$, $u$ is the inverse Fourier transform of the function $\dfrac{(\FF f)(\xi)}{|\xi|_p^\alpha +\mu}$ belonging to $L^2(\Qp)$. The function $\DA \psi +\mu \psi$ is the inverse Fourier transform of the function $\left( |\xi|_p^\alpha +\mu\right) (\FF \psi )(\xi)$, also belonging to $L^2(\Qp)$. By the Plancherel formula,
$$
\langle u,\DA \psi +\mu \psi \rangle =\int\limits_{\Qp}(\FF f)(\xi )\overline{\FF \psi)(\xi )}\,d\xi
$$
where $\FF f$ is bounded, $\FF \psi \in \mathcal D(\Qp)$.

In particular, $(\FF \psi )(\xi )=0$ for $|\xi|_p\ge p^N$, with some $N\in \mathbb N$. Then for any $n\ge N$
\begin{multline*}
\langle u,\DA \psi +\mu \psi \rangle =\int\limits_{\Qp}f(x)\,dx\int\limits_{\Qp}\overline{\psi (y)}\,dy\int\limits_{|\xi |_p\le p^n}\chi ((x-y)\xi )\,d\xi \\
=p^n\iint\limits_{|x-y|_p\le p^{-n}}f(x)\overline{\psi (y)}\,dx\,dy=p^n\int\limits_{\Qp}f(x)\,dx\int\limits_{|z|_p\le p^{-n}}\overline{\psi (x-z)}\,dz.
\end{multline*}

Since $\psi$ is locally constant, and its locally constancy exponent can be chosen the same for a neighborhood of every point $x$, for a sufficiently large $n$ we have
$$
p^n \int\limits_{|z|_p\le p^{-n}}\overline{\psi (x-z)}\,dz=\overline{\psi (x)},
$$
so that $\langle u,\DA \psi +\mu \psi \rangle =\langle f,\psi \rangle$.

Thus, for every $u\in D(A)$, $Au=\DA u$ in the distribution sense. Suppose that $u\in \LL$, $\DA u\in \LL$  in the distribution sense. Set $f=(\DA +\mu I)u$ (in the distribution sense). Denote $u'=R_\mu (A)f$. Then $u'\in D(A)$ and $(\DA +\mu )(u-u')=0$ in the distribution sense. $\DA$ is a convolution operator, and the existence of the convolution of distributions from $\mathcal D'(\Qp)$ implies \cite{VVZ} the following equality for their Fourier transforms:
$$
\left( |\xi|_p^\alpha +\mu\right) ((\FF u)(\xi )-(\FF u')(\xi))=0
$$
for all $\xi \in \Qp$. Therefore $u=u'\in D(A)$. $\qquad\blacksquare$.

\bigskip
\section{$L^1$-Theory of the Vladimirov Type Operator on a $p$-Adic Ball}

{\it 4.1. The Heat-Like Semigroup}.

On a ball $B_N$, $N\in \mathbb Z$, we consider the Cauchy problem (\ref{2.9})-(\ref{2.10}). Its fundamental solution (\ref{2.11}) defines a contraction semigroup
$$
(T_N(t)u)(x)=\int\limits_{B_N}Z_N(t,x-\xi )u(\xi )\,d\xi
$$
on $L^1(B_N)$.

\medskip
\begin{prop}
The semigroup $T_N$ is strongly continuous.
\end{prop}

\medskip
{\it Proof}. For $u\in L^1(B_N)$, we have
$$
\left\| T_N(t)u-u\right\|_{L^1(B_N)}\le I_1(t)+I_2(t)
$$
where
$$
I_1(t)=\int\limits_{B_N}dx\left| \int\limits_{B_N}e^{\lambda t}Z_N(t,x-\xi )u(\xi )\,d\xi -u(x)\right| ,
$$
$$
I_2(t)=p^N|c(t)|\int\limits_{B_N}|u(\xi )|\,d\xi .
$$
Since $c(t)\to 0$, as $t\to 0$, we get $I_2(t)\to 0$.

For small values of $t$, we write
\begin{multline*}
I_1(t)=\int\limits_{B_N}dx\left| \int\limits_{B_N}Z(t,x-\xi )u(\xi )\,d\xi -u(x)+\int\limits_{B_N}(e^{\lambda t}-1)Z(t,x-\xi )u(\xi )\,d\xi \right| \\
\le \int\limits_{B_N}dx\left| \int\limits_{B_N}Z(t,x-\xi )u(\xi )\,d\xi -u(x)\right| +Ct\int\limits_{B_N}dx \int\limits_{B_N}Z(t,x-\xi )|u(\xi )|\,d\xi \\
\stackrel{\text{def}}{=}I_{1,1}(t)+I_{1,2}(t).
\end{multline*}

By the Young inequality and the identity (\ref{3.1}), extending $u$ by zero to a function $\widetilde{u}$ on $\Qp$ we obtain
$$
I_{1,2}(t)\le Ct\int\limits_{\Qp}dx \int\limits_{\Qp}Z(t,x-\xi )|\widetilde{u}(\xi )|\,d\xi \le Ct\|\widetilde{u}\|_{\LL}=Ct\|u\|_{L^1(B_N)}\to 0,
$$
as $t\to 0$.

Next,
\begin{multline*}
I_{1,1}(t)=\int\limits_{B_N}dx\left| \int\limits_{\Qp}Z(t,x-\xi )\widetilde{u}(\xi )\,d\xi -\widetilde{u}(x)\right| \le \int\limits_{\Qp}dx\left| \int\limits_{\Qp}Z(t,x-\xi )\widetilde{u}(\xi )\,d\xi -\widetilde{u}(x)\right| \\
=\|S(t)\widetilde{u}-\widetilde{u}\|_{\LL}\to 0,
\end{multline*}
as $t\to 0$, by the $C_0$-property of $S(t)$. $\qquad \blacksquare$

\bigskip
{\it 4.2. The Generator}.

Denote by $A_N$ the generator of the contraction semigroup $T_N$ on $L^1(B_N)$. By the Hille-Yosida theorem, $A_N$ has a bounded resolvent $(A_N+\mu I)^{-1}$ for each $\mu >0$. In order to study the domain $D(A_N)$, we need the following auxiliary result.

\medskip
\begin{lem}
Let the support of a function $u\in \LL$ be contained in $\Qp \setminus B_N$. Then the restriction to $B_N$ of the distribution $\DA u\in \DDD$ coincides with the constant
$$
R_N=R_N(u)=\frac{1-p^\alpha}{1-p^{-\alpha -1}}\int\limits_{|x|_p>p^N}|x|_p^{-\alpha -1}u(x)\,dx.
$$
\end{lem}

\medskip
{\it Proof}. Let $\psi\in \mathcal D(B_N)$. Then $\langle \DA u,\psi\rangle =\langle u,\DA \psi \rangle$ where
$$
\left( \DA \psi \right) (x)=\frac{1-p^\alpha }{1-p^{-\alpha -1}}\int\limits_{\Qp}|y|_p^{-\alpha -1}[\psi (x-y)-\psi (x)]\,dy.
$$

Since $\psi (x)=0$ for $|x|_p>p^N$, we find using the ultrametric property that
\begin{multline*}
\langle u,\DA \psi \rangle =\frac{1-p^\alpha }{1-p^{-\alpha -1}}\int\limits_{|x|_p>p^N}u(x)\,dx\int\limits_{\Qp}|y|_p^{-\alpha -1}[\psi (x-y)-\psi (x)]\,dy \\
= \frac{1-p^\alpha }{1-p^{-\alpha -1}}\int\limits_{|x|_p>p^N}u(x)\,dx\int\limits_{|z|_p\le p^N}|x-z|_p^{-\alpha -1}\psi (z)\,dz\\
=\frac{1-p^\alpha }{1-p^{-\alpha -1}}\int\limits_{|x|_p>p^N}|x|_p^{-\alpha -1}u(x)\,dx\int\limits_{|z|_p\le p^N}\psi (z)\,dz.
\end{multline*}
Since $\psi \in \mathcal D(B_N)$ is arbitrary, this implies the required property. $\qquad \blacksquare$

\medskip
Now we can formulate the main result of this section. As before, $A$ denotes the generator of the semigroup $S(t)$ on $\LL$.
\begin{prop}
If $\psi \in D(A)$, then the restriction $\psi_N$ of the function $\psi$ to $B_N$ belongs to $D(A_N)$, and $A_N\psi_N=\left( \DAN -\lambda \right)\psi_N$ where $\DAN \psi_N$ is understood in the sense of $\mathcal D'(B_N)$, that is $\psi_N$ is extended by zero to a function on $\Qp$, $\DA$ is applied to it in the distribution sense, and the resulting distribution is restricted to $B_N$.
\end{prop}

\medskip
{\it Proof.} For $\psi \in D(A)$, we have to check that $\DAN \psi_N \in L^1(B_N)$ and
\begin{equation}
\label{4.1}
\left\| -\frac1{t}[T_N(t)\psi_N-\psi_N]-\left( \DAN -\lambda \right)\psi_N\right\|_{L^1(B_N)}\longrightarrow 0,
\end{equation}
as $t\to +0$.

By Lemma 5, we can write on $B_N$ as follows: $\psi =\psi_N+(\psi -\psi_N)$, $A\psi=\DAN \psi_N+R_N$, $R_N=R_N(\psi -\psi_N)$, so that $\DAN \psi_N=A\psi -R_N\in L^1(B_N)$.

Next, it follows from (\ref{2.11}) that
\begin{multline}
\label{4.2}
(T_N(t)\psi_N)(x)=\int\limits_{B_N}Z(t,x-y)\psi (y)\,dy +c(t)\int\limits_{B_N}\psi (y)\,dy+\lambda t \int\limits_{B_N}Z(t,x-y)\psi (y)\,dy\\
+d(t)\int\limits_{B_N}Z(t,x-y)\psi (y)\,dy
\end{multline}
where $d(t)=O(t^2)$, $t\to 0$. By the strong continuity of $S(t)$,
$$
\left\| -\lambda \int\limits_{B_N}Z(t,x-y)\psi (y)\,dy+\lambda \psi_N(x)\right\|_{L^1(B_N)}\longrightarrow 0,
$$
as $t\to 0$. It follows from the Young inequality that
$$
\frac1{t}\left\| d(t)\int\limits_{B_N}Z(t,x-y)\psi (y)\,dy\right\|_{L^1(B_N)}\longrightarrow 0.
$$

Next, it is checked directly that $c(0)=c'(0)=0$, so that $c(t)=O(t^2)$, $t\to 0$, and the contribution of the second term in (\ref{4.2}) is negligible. Thus it remains to consider the first term, that is
$$
V(t,x)=\int\limits_{B_N}Z(t,x-y)\psi (y)\,dy=v_1(t,x)-v_2(t,x),\quad x\in B_N,
$$
where
$$
v_1(t,x)=\int\limits_{\Qp}Z(t,x-y)\psi (y)\,dy,
$$
$$
v_2(t,x)=\int\limits_{|y|_p>p^N}Z(t,x-y)\psi (y)\,dy.
$$

Using the representation (\ref{2.7}) and noticing that $|x-y|_p=|y|_p$, if $|x|_p\le p^N$, $|y|_p>p^N$, we find that
$$
v_2(t,x)=\sum\limits_{m=1}^\infty \frac{(-1)^m}{m!}\cdot \frac{1-
p^{\alpha m}}{1-p^{-\alpha m-1}}t^m\int\limits_{|y|_p>p^N}|y|_p^{-\alpha m-1}\psi (y)\,dy
$$
($v_2=v_2(t)$ does not depend on $x\in B_N$). In particular, as $t\to 0$,
\begin{equation}
\label{4.3}
-\frac1{t}v_2(t)\longrightarrow -\frac{1-p^{\alpha }}{1-p^{-\alpha -1}}\int\limits_{|y|_p>p^N}|y|_p^{-\alpha m-1}\psi (y)\,dy
\end{equation}
where the convergence can be interpreted as the one in $L^1(B_N)$. The right-hand side of (\ref{4.3}) coincides with $-R_N$, and since
\begin{equation}
\label{4.4}
\DAN \psi_N=A\psi -R_N,
\end{equation}
we obtain that
\begin{multline*}
\left\| -\frac1{t}[T_N(t)\psi_N-\psi_N]-\left( \DAN -\lambda \right)\psi_N\right\|_{L^1(B_N)}\le \left\| -\frac1{t}[S(t)\psi -\psi ]-A\psi\right\|_{L^1(B_N)}+o(1)\\
\le \left\| -\frac1{t}[S(t)\psi -\psi ]-A\psi\right\|_{L^1(\Qp)}+o(1)\longrightarrow 0,
\end{multline*}
as $t\to 0$, which completes the proof. $\qquad \blacksquare$

\bigskip
\section{Nonlinear Equations: the Main Result}

Let us return to Eq. (\ref{1.1}) interpreted as Eq. (\ref{1.2}) on $\LL$, where the linear operator $A$ is a generator of the semigroup $S(t)$. As before, $\varphi$ is a strictly monotone increasing continuous real function, $|\varphi (s)|\le C|s|^m$, $m\ge 1$. Below we re-interpret Eq. (\ref{1.1}) as the equation
\begin{equation}
\label{5.1}
\frac{\partial u}{\partial t}+\overline{A\varphi}(u)=0
\end{equation}
where $\overline{A\varphi}$ is the closure of $A\varphi$. It follows from Proposition 2 that already the operator $A\varphi$ is densely defined, the more so it is valid for $\overline{A\varphi}$.

Recall that a mild solution of the Cauchy problem for a nonlinear equation with the initial condition $u(0,x)=u_0(x)$ is defined as a function given by a limit, uniformly on compact time intervals, of solutions of the problem for the difference equations approximating the differential one. This the usual ``nonlinear version'' of the notion of a generalized solution; see \cite{Ba} for the details.

\medskip
\begin{teo}
The operator $\overline{A\varphi}$ is $m$-accretive, so that, for any initial function $u_0\in \LL$, the Cauchy problem for Eq. (\ref{5.1}) has a unique mild solution.
\end{teo}

\medskip
{\it Proof}. The second statement is a consequence of the first one and the Crandall-Liggett theorem; see Theorem 4.3 in \cite{Ba}.

The accretivity of $A\varphi$ means that the inequality
$$
\|x-y\|_{\LL}\le \|(I+A\varphi )(x)-(I+A\varphi )(y)\|_{\LL}
$$
is valid on its domain. Therefore, in order to prove the $m$-accretivity of $\overline{A\varphi}$, it suffices to check that the operator $I+A\varphi$ has a dense range. Since $\DD$ is dense in $\LL$, it will be sufficient to prove a little more, namely the solvability of Eq. (\ref{2.15}) for any $f\in \LL \cap L^\infty (\Qp)$.

For such a function $f$, we consider Eq. (\ref{2.15}) and the regularized equation Eq. (\ref{2.16}). Using Eq. (\ref{2.19}) we find that
\begin{equation}
\label{5.2}
|v_\varepsilon (x)|\le \beta^{-1}\left( \|f\|_{L^\infty (\Qp)}\right)
\end{equation}
for almost all $x\in \Qp$. This means that for any fixed $N$,
$$
|R_N(v_\varepsilon )|\le C
$$
where $C$ does not depend on $\varepsilon$, so that the set of constant functions $\{ R_N(v_\varepsilon ),0<\varepsilon<1\}$ is relatively compact in $L^1(B_N)$.

On the other hand, it follows from (\ref{2.17}), (\ref{2.18}) (with $\Omega=\Qp$) and the translation invariance of $A$ that the family of functions $w_\varepsilon =f-Av_\varepsilon$ satisfies the inequalities
\begin{equation}
\label{5.3}
\|w_\varepsilon\|_{\LL}\le \|f\|_{\LL};
\end{equation}
\begin{equation}
\label{5.4}
\int\limits_{\Qp}|w_\varepsilon (x+h)-w_\varepsilon (x)|\,dx \le \int\limits_{\Qp}|f(x+h)-f(x)|\,dx
\end{equation}
(for any $h\in \Qp$).

The conditions (\ref{5.3}) and (\ref{5.4}) imply the relative compactness of $\{ w_\varepsilon\}$, thus of $\{ Av_\varepsilon\}$, in $L^1_{\text{loc}}(\Qp)$, that is the compactness of the closure of the restriction $(Av_\varepsilon)\upharpoonright_X$ for any bounded measurable subset $X\subset \Qp$. This is a consequence of the criterion for relative compactness in $L^1(G)$ where $G$ is a compact group (see Theorem 4.20.1 in \cite{Ed}) applied to the case $G=B_N$ (the additive group of a $p$-adic ball).

Denote by $v_{\varepsilon,N}$ the restriction of $v_\varepsilon$ to $B_N$. It follows from (\ref{4.4}) that the set $\left\{ \DAN v_{\varepsilon,N}\right\}$ is relatively compact in $L^1(B_N)$. Since $\DAN =A_N+\lambda$ has a bounded inverse on $L^1(B_N)$, this implies the relative compactness of $\{ v_{\varepsilon,N}\}$ in $L^1(B_N)$ for each $N$, thus the one of $\{v_\varepsilon\}$ in $L^1_{\text{loc}}(\Qp)$. Let $v$ be a limit point. Together with the relative compactness of $\{ Av_\varepsilon\}$, the above reasoning proves the solvability of Eq. (\ref{2.15}) because, by Fatou's lemma and (\ref{5.2}), a limit point of $\{ Av_\varepsilon\}$ belongs to $\LL$. Therefore $\beta (v)\in \LL$. By (\ref{5.2}), $v\in L^\infty (\Qp)$, so that $\beta (v)\in L^\infty (\Qp)$, $v=\varphi (\beta (v))$, $|v(x)|\le C|\beta (v)|^m\le C_1|\beta (v)|$, and $v$ belongs to $\LL$. $\qquad \blacksquare$

\bigskip
\section{Explicit Solution: an Example}

Let us consider Eq. (\ref{1.1}) with $\alpha >0$, $\varphi (u)=|u|^m$, $m>1$. We look for a solution of the form
\begin{equation}
\label{6.1}
u(t,x)=\rho \left( \frac{|x|_p^\gamma}{t_0-t}\right)^\nu,\quad 0<t<t_0, x\in \Qp,
\end{equation}
where $t_0>0$, $\gamma >0$, $\nu >0$, $0\ne \rho \in \mathbb R$.

We have
$$
\frac{\partial u}{\partial t}=\nu \rho |x|_p^{\gamma \nu}(t_0-t)^{-\nu -1},
$$
$$
\DA (|u|^m)=|\rho |^m (t_0-t)^{-\nu m}\DA \left( |x|_p^{\gamma \nu m}\right).
$$
Comparing powers of $t_0-t$ we find that $\nu =(m-1)^{-1}$.

We understand $\DA$ as the convolution operator (\ref{3.13}) defined, in particular, for any $\alpha >0$ by analytic continuation; $f_\beta$ is defined for all real $\beta$ except $\beta =1$. If $\beta \ne 1$, $-\alpha +\beta \ne 1$, then $f_{-\alpha }*f_\beta =f_{-\alpha +\beta}$ where the convolution is understood in the sense of $\DDD$. See \cite{AKS,K2001,VVZ} for the details.

The above convolution identity can be written as
$$
\DA (|\cdot |_p^{\beta -1})=\Gamma_p(\beta )\frac{|\cdot |_p^{-\alpha +\beta -1}}{\Gamma_p(\beta -\alpha)}
$$
or, if we substitute $\beta +1$ for $\beta$, as
\begin{equation}
\label{6.2}
\DA (|\cdot |_p^\beta )=\frac{\Gamma_p(\beta +1)}{\Gamma_p(\beta -\alpha +1)}|\cdot |_p^{\beta -\alpha},\quad \beta \ne \alpha .
\end{equation}

Calculating $\DA (|u|^m)$ by (\ref{6.2}), substituting into (\ref{6.1}) and canceling powers of $t_0-t$ we come to the identity
$$
\frac{\rho}{m-1}|x|_p^{\frac{\gamma}{m-1}}=-|\rho |^m\frac{\Gamma_p(\frac{\gamma m}{m-1}+1)}{\Gamma_p(\frac{\gamma m}{m-1}-\alpha +1)}|x|_p^{\frac{\gamma m}{m-1}-\alpha }
$$
implying $\gamma =\alpha$,
\begin{equation}
\label{6.3}
\frac{\rho}{m-1}=-|\rho |^m\frac{\Gamma_p(\frac{\alpha m}{m-1}+1)}{\Gamma_p(\frac{\alpha }{m-1}+1)}.
\end{equation}

Both the numerator and denominator in the right-hand side of (\ref{6.3}) are negative. Thus $\rho <0$, and we obtain our solution in the form
\begin{equation}
\label{6.4}
u(t,x)=\rho (t_0-t)^{-\frac1{m-1}}|x|_p^{\frac{\alpha}{m-1}}
\end{equation}
where
$$
\rho =-\left[ \frac{\Gamma_p(1+\frac{\alpha}{m-1})}{(m-1)\Gamma_p(1+\frac{\alpha m}{m-1})}\right]^{\frac1{m-1}}.
$$

In a similar way, we can obtain instead of the solution (\ref{6.4}), the solution
$$
u(t,x)=\mu (t_0+t)^{-\frac1{m-1}}|x|_p^{\frac{\alpha}{m-1}},\quad t>0,x\in \Qp,
$$
where $\mu =-\rho$.

\section*{Acknowledgments}
The second author is grateful to Mathematical Institute, Linnaeus University, for hospitality during his visits to V\"axj\"o. The work of the second author was also supported in part by Grant 23/16-18 ``Statistical dynamics, generalized Fokker-Planck equations, and their applications in the theory of complex systems'' of the Ministry of Education and Science of Ukraine.

\end{document}